\documentclass[leqno,twoside,a4paper,12pt]{article}
\usepackage{amsmath,amssymb,amsthm}
\usepackage{url}
\usepackage[midshaft,scriptlabels]{diagrams}

\DeclareFontFamily{U}{UWCyr}{}
\DeclareFontShape{U}{UWCyr}{m}{n}{%
  <5> <6> <7> <8> <9>
  <10> <10.95> <12> <14.4> <17.28> <20.74> <24.88> wncyr10
  }{}
\DeclareMathAlphabet{\cyrm}{U}{UWCyr}{m}{n}
\DeclareSymbolFont{cyrm}{U}{UWCyr}{m}{n}
\DeclareSymbolFontAlphabet{\cyrm}{cyrm}
\DeclareMathSymbol{\Evo}{\cyrm}{cyrm}{"03}
\newtheorem{theorem}{Theorem}
\newtheorem{lemma}{Lemma}
\newtheorem{corollary}{Corollary}
\newtheorem{proposition}{Proposition}

\theoremstyle{definition}
\newtheorem{definition}{Definition}

\theoremstyle{remark}
\newtheorem{remark}{Remark}

\DeclareMathOperator{\byd}{\raisebox{-.2ex}{$\overset{\text{\tiny
def}}{=}$}}

\DeclareMathOperator{\im}{im}

\DeclareMathOperator{\pro}{pr}

\DeclareMathOperator{\id}{id}

\DeclareMathOperator{\Imm}{Imm}

\DeclareMathOperator{\Sym}{Sym}

\DeclareMathOperator{\ext}{ext}

\DeclareMathOperator{\inter}{int}

\DeclareMathOperator{\Aff}{Aff}

\newarrow{TeXto}----{->}
\newarrow{MapsTo}|--->
\newcommand{\mto}{\rMapsTo[]\relax}
\newarrow{hTeXto}C---{->}
\newarrow{hTo}C--->

\newcommand{\cprime}{\/{\mathsurround=0pt$'$}}

\newcommand{\bsi}{\boldsymbol{\sigma}}

\newcommand{\R}{\mathbb{R}}

\newcommand{\cC}{\mathcal{C}}

\newcommand{\cE}{\mathcal{E}}

\newcommand{\cV}{\mathcal{V}}
\newcommand{\cH}{\mathcal{H}}

\newcommand{\cS}{\mathcal{S}}

\newcommand{\cW}{\mathcal{W}}

\newcommand{\eval}[2][\right]{\relax
  \ifx#1\right\relax \left.\fi#2#1\rvert}
\title{On the geometry of Grassmannian equivalent connections}
\author{Gianni Manno
\\ {\footnotesize Dip. di Matematica ``Ennio De Giorgi'', Universit\`a di Lecce}
\\ {\footnotesize email: Gianni.Manno@unile.it}}

\date{}
\begin{document}

\maketitle

\begin{abstract}

We introduce the equation of $n$-dimensional totally geodesic
submanifolds of $E$ as a submanifold of the second order jet space
of $n$-dimensional submanifolds of $E$. Next we study the geometry
of $n$-Grassmannian equivalent connections, that is linear
connections without torsion admitting the same equation of
$n$-dimensional totally geodesic submanifolds. We define the
$n$-Grassmannian structure as the equivalence class of such
connections, recovering for $n=1$ the case of theory of
projectively equivalent connections. By introducing the equation
of parametrized $n$-dimensional totally geodesic submanifolds as a
submanifold of the second order jet space of the trivial bundle on
the space of parameters, we discover a relation of covering, in
the sense of \cite{Many99,KrVe98,KrVi84,KrVi89}, between the
`parametrized' equation and the `unparametrized' one. After having
studied symmetries of these equations, we discuss the case in
which the space of parameters is equal to $\R^n$.

\end{abstract}

%\tableofcontents

\section{Introduction}

In the 1920's Weyl introduced in \cite{Weyl21} the notion of
projectively equivalent connections, that is torsion-free linear
connections on a manifold $E$ having the same geodesics as paths,
namely up to reparametrization. Locally, two such connections
$\Gamma$ and $\Gamma^{\prime}$ are related by
\begin{equation}\label{eq.proj.equiv.connecions.locally}
\Gamma_A{}^C{}_B-{\Gamma'}_A{}^C{}_B=\delta^C_A\Phi_B+\delta^C_B\Phi_A.
\end{equation}
where $\Phi_A,\Phi_B\in C^\infty(E)$.  This approach was pursued
by Eisenhart \cite{Eisen22}, Veblen \cite{Veblen22} and Thomas
\cite{Thomas25,Thomas26}. He introduced, on an $l$-dimensional
manifold with connection $\Gamma$, the projective invariants
\begin{equation}\label{eq.projective.invariants}
\Pi_A{}^C{}_B=\Gamma_A{}^C{}_B-\frac{1}{l+1}
\left(\Gamma_A{}^F{}_F \, \delta^C_B+\Gamma_B{}^F{}_F \,
\delta^C_A\right),
\end{equation}
that is quantities which do not change for a connection which is
projectively equivalent to $\Gamma$. He called them collectively a
\emph{projective connection}.

In opposition with these mathematicians, Cartan in \cite{Cartan24}
approached the concept of projective connection in a different
way. He imagined to attach to every point of a manifold a
projective space of the same dimension, and he gave a rule to
connect the projective spaces corresponding to infinitesimally
near points. He defined the geodesics as curves whose
``developments'' are straight lines in the corresponding
projective spaces. Then in \cite{KobNagano} the authors gave a
rigorous foundation of Cartan's approach and explained the
relationships with the point of view of Thomas. They introduced
the \emph{projective structure} on a manifold $E$ as a particular
subbundle of the second frame bundle of $E$, and a
\emph{projective symmetry} as a diffeomorphism of $E$ preserving
the projective structure. Then they regarded a projective
connection as a particular connection on the projective structure.

%\smallskip
In this paper we approach the theory of projectively equivalent
connections by using jets of submanifolds (also known as manifold
of higher contact elements)
\cite{MaThesis,MaAnkara,MaVi02,MaVi03,Vin81,Vin84,Vin88}. Roughly
speaking, the $r$-jet space of $n$-dimensional submanifolds
$J^r(E,n)$ of a given manifold $E$ is the space of equivalence
classes of $n$-dimensional submanifolds of $E$ having an $r$-th
contact at a certain point. For instance, $J^0(E,n)=E$ and
$J^1(E,n)$ is the Grassmann bundle on $E$ of $n$-dimensional
subspaces of $TE$. The advantage of using jets of submanifolds is
that they provide a differential calculus on unparametrized
submanifolds, so that any object constructed on them, like
connections, vector fields, etc..., are independent of the Lie
group $G^r_n$ of $r$-jets of `reparametrizations'
\cite{Ehr52,GrKr98,KMS93,Kru01}. This makes jets of submanifolds
natural geometrical objects for studying differential properties
of the geometry of submanifolds where the parametrization is not
important (like projective geometry of connections). Then, in our
scheme, we do not need to use principal bundles like frame
bundles, as in \cite{KobNagano}, to develop an intrinsic theory of
projectively equivalent connections, obtained here in section
\ref{sec.Totally.Geod.Subm.in.the.Jet.of.Subm} as the particular
case $n=1$. Taking this into account, we are able to characterize
the projective structure on $E$ associated with a torsion-free
linear connection $\Gamma$ on $TE\rTo TE$ as a particular
distribution on the projectivized tangent bundle $J^1(E,1)$ of
$E$, or, equivalently, as a particular section
$\ddot{\Gamma}\colon J^1(E,1)\rTo J^2(E,1)$.

Our geometrical setting allows us to generalize in a natural
manner our reasoning to the case of $n$-dimensional submanifolds.
In section \ref{sec.Totally.Geod.Subm.in.the.Jet.of.Subm} we
introduce the notion of \emph{$n$-Grassmannian equivalent
connections} associated with a torsione-free linear connection as
the equivalence class of connections admitting the same equation
of  unparametrized $n$-dimensional totally geodesic submanifolds.
We find a distinguished representative $\ddot{\Gamma}$ of such a
class as a section $\ddot{\Gamma}\colon J^1(E,n)\rTo J^2(E,n)$, or
equivalently as a particular distribution on the Grassmann bundle
$J^1(E,n)$, that we call an \emph{$n$-Grassmannian structure},
recovering for $n=1$ the case of projectively equivalent
connections. We show that the image of $\ddot{\Gamma}$ coincides
with the equation $\cS$ of unparametrized totally geodesic
submanifolds, so that we treat the Grassmannian structure, the
class of Grassmannian equivalent connections and the equation of
unparametrized totally geodesic submanifolds practically as the
same object. In this way we obtain a more simple description of
this geometry than that given in \cite{Dhooghe94}, where the
author generalizes definitions and results of \cite{KobNagano} to
the Grassmannian case by straightforward computations with the
same approach therein. Also, taking into consideration the
particular structure of $\cS$, we calculate the $n$-Grassmannian
invariants, that is quantities which do not change for
$n$-Grassmannian equivalent connections, and we see as in the case
$n=1$ such invariants reduce to those introduced by Thomas.

We construct the map $\ddot{\Gamma}$ through a connection
$\dot{\Gamma}$ on $J^1(E,n)$ which we prove is naturally
associated with $\Gamma$. We take inspiration from
\cite{JanyModyGeneral}, where such a construction is performed in
the case of jets of time-like curves of an oriented manifold. Our
result is much more general, as we construct $\dot{\Gamma}$ in the
case of $n$-dimensional submanifolds which are not required to be
oriented. Moreover, in section \ref{sec.covering}, we find another
way to get such a connection.

In section \ref{sec.Totally.Geod.Subm.in.the.Jet.Bundle} we define
the equation $\cS_{\pro_{M}}$ of parametrized totally geodesic
submanifolds as a submanifold of $J^2{\pro_{M}}$ where
${\pro_{M}}\colon M\times E\rTo M$, and M is the space of
`parameters'. Following the reasoning of section
\ref{sec.Totally.Geod.Subm.in.the.Jet.of.Subm}, we define a
connection $\dot{\Gamma}_{\pro_{M}}$ on $J^1{\pro_{M}}\rTo M\times
E$ and an operator $\ddot{\Gamma}_{\pro_{M}}$.

In section \ref{sec.symm.of.geodesic.equation}, we discuss
symmetries of both the equations $\cS$ and $\cS_{\pro_{M}}$. We
see that the most general infinitesimal symmetry of the equation
$\cS$ is a vector field on $J^1(E,n)$ which preserves the
$n$-Grassmannian structure $\ddot{\Gamma}$, that we call a
\emph{contact Grassmannian symmetry}. If we restrict our attention
to the case $n=1$, we realize that this point of view is more
general that the point of view of standard projective symmetries
studied in literature (see \cite{Aminova01,Aminova04} and
references therein), where just vector fields on $E$ have been
considered.

In the case $n=1$, in \cite{Aminova01} the author find, starting
from the equation of ``parametrized'' geodesics, the equation of
unparametrized geodesics by a straightforward substitution. In
section \ref{sec.covering}, we prove that such a substitution is a
local description of a natural geometrical property of the
equation of unparametrized totally geodesic submanifolds which is
that of possessing a \emph{covering equation}, in the sense of
\cite{Many99,KrVe98,KrVi84,KrVi89}. More precisely, we find a
natural surjection from $\cS_{\pro_{M}}$ to $\cS$ which we prove
preserves the differential structures of the equations. Also, we
show that we can obtain the connection $\dot{\Gamma}$ as the
quotient connection of $\dot{\Gamma}_{\pro_{M}}$ via this covering
map. Finally, in section \ref{sec.caso.particolare} we prove that
in the case $M=\R^n$ this covering can be obtained by factoring
the equation $\cS_{\pro_{\R^n}}$ by the affine group of $\R^n$,
which turns to be a subgroup of the group of symmetries of
$\cS_{\pro_{\R^n}}$.

\section{Geometry of Differential Equations on Submanifolds}\label{sec.finite.jet.spaces}
\label{sec.preliminary.definitions}

In this sections we recall basic notions of the theory of jets of
submanifolds, differential equations and their symmetries. Our
main sources are
\cite{Many99,KrVe98,MaThesis,MaAnkara,MaVi02,MaVi03,MoVi94,Olver01,Vin81,Vin84,Vin88,VinII}.

\subsection{Jet Spaces}
Let $E$ be an $(n+m)$-dimensional smooth manifold and $L$ an
$n$-dimensional embedded submanifold of $E$. Let $(u^A)$ be a
local chart on $E$. The chart $(u^A)$ can be divided in two parts,
$(u^A)=(u^\lambda,u^i)$, $\lambda=1\dots n$ and $i=n+1\dots m+n$,
such that the submanifold $L$ is locally described by
$u^i=f^i(u^{1},u^{2},\ldots ,u^{n})$.

The chart $(u^\lambda,u^i)$ is said to be a \emph{divided chart}
which is \emph{concordant} to $L$. Here, and in the rest of the
paper, Greek indices run from $1$ to $n$, Latin indices run from
$n+1$ to $m+n$, and capital letters from $1$ to $m+n$. Otherwise
we shall specify them. Also, all submanifolds are \emph{embedded}
submanifolds.

Let $L'$ be another $n$-dimensional submanifolds locally described
$u^i=f'^i(u^{1},u^{2},\ldots ,u^{n})$. We say that $L$ and $L'$
have a \emph{contact of order $r$} at $p$ if $f$ and $f'$ have a
contact of order $r$ at $p$. Locally this means that the Taylor
expansion of $(f-f')$ around $p$ vanishes up to order $r$. This
property is invariant by coordinate transformations.

The above relation is an equivalence relation; an equivalence
class is denoted by~$[L]^{r}_{p}$. The set of such classes is said
to be the \emph{$r$-jet of $n$-dimensional submanifolds of $E$}
and it is denoted by $J^r(E,n)$.

The set $J^r(E,n)$ has a natural manifold structure. Namely, let
$\bsi =(\sigma _{1},\sigma _{2},\ldots,\sigma_{k})$, with
$1\leq\sigma_i\leq n$ and $r\in\mathbb{N}$, be a multi-index, and
$\left| \bsi \right|\byd k$. Any divided chart $(u^\lambda,u^i)$
at $p\in E$ induces the local chart $\left(u^{\lambda
},u_{\bsi}^{i}\right)$ at $[L]_{p}^{r}\in J^{r}(E,n)$, where
$|\bsi|\leq r$ and the functions $u^j_{\bsi}$ are determined by
$u_{\bsi}^{i}\circ j_{r}L=\partial^{|\bsi|}f^{i}\big/\partial
u^{\bsi}$.

We note that $J^{0}(E,n)=E$. Also, $[L]_{p}^{1}$ can be identified
with the tangent space $T_{p}L$. So $J_{p}^{1}(E,n)$ is the
Grassmannian of $n$-dimensional subspaces of $T_{p}E$, and
$J^{1}(E,n)$ is the Grassmann bundle on $E$ of $n$-dimensional
subspaces of $TE$.

\smallskip
We have the following natural maps:
\begin{enumerate}
\item the embedding $j_{r}L\colon L \rTo J^{r}(E,n),\,\, p \mto
[L]_{p}^{r}$\,,

\item the projection $\pi_{k,h}\colon J^{k}(E,n)\rTo J^{h}(E,n),\,[L]^k_p\mto
[L]^h_p\,\quad k\geq h$.
\end{enumerate}

We denote by $L^{(r)}$ the image of $j_rL$. We call the tangent
plane $T_{[L]^r_p}L^{(r)}$ an \emph{$R$-plane}.

\smallskip
We denote by $\chi(J^r(E,n))$ the module of vector fields on
$J^r(E,n)$.

\smallskip
In the rest of the paper we shall put $\theta_r=[L]^r_p$.

\smallskip
The Cartan plane $\cC^r_{\theta_r}$ on $J^r(E,n)$ at $\theta_r$ is
defined as the span of the planes $T_{\theta _{r}}L^{(r)}$ with
$L$ running over all $n$-dimensional submanifolds of $E$. The
correspondence $\theta_r\mto \cC^r_{\theta_r}$ is called the
\emph{Cartan distribution}. We denote by $\cC^r(D)$ the set of
vector fields lying in the Cartan distribution of $J^r(E,n)$. It
is easy to realize that to each point $\theta_{r+1}$ there
corresponds the $R$-plane $R_{\theta_{r+1}}=T_{\theta_r}L^{(r)}$,
and that
\begin{equation}\label{eq.splitting.Cartan}
\cC^r_{\theta_r}=R_{\theta_{r+1}}\oplus \ker
T_{\theta_r}\pi_{r,r-1}.
\end{equation}
This direct sum is not canonical as we have many $R$-planes
passing through a point $\theta_r$. A diffeomorphism of $J^r(E,n)$
which preserves the Cartan distribution is called a \emph{contact
transformation}. Analogously we can define \emph{contact fields}
on $J^r(E,n)$. We denote by $D_{\cC^r}$ the set of such fields.

A \emph{differential equation $\mathcal{E}$ of order $r$ on
$n$-dimensional submanifolds of a manifold E} is a closed
submanifold of $J^r(E,n)$. A \emph{solution} is a $n$-dimensional
submanifold $L$ of $E$ such that $L^{(r)}\subset\cE$.

Let $\eta\colon F\rTo J^1(E,n)$ be a vector bundle and let
$\eta_r=\pi^*_{r,1}(\eta)$ be the pull-back bundle of $\eta$
through the map $\pi^*_{r,1}$. A differential equation $\cE\subset
J^r(E,n)$ can be described by $\phi=0$ where
$\phi\in\Gamma(\eta_r)$. We can associate with $\phi$ the
following operation:
\begin{equation}\label{rem.jet.morphism1}
\square_{\phi}\colon L\mto (\eta^*(\pi_{r,1})\circ\phi\circ
j_rL)\in\Gamma(\eta|_L),\,\,\,L\subset E\,,
\end{equation}
that we call a \emph{non-linear differential operator} following
the terminology of \cite{VinII}. Conversely, any map
$\square\colon L\mto \square(L)\in\Gamma(\eta|_L)$ such that
$\square(L)(p)=\square(\tilde{L})(p)$ whenever
$[L]^r_p=[\tilde{L}]^r_p$, induces the following section:
\begin{equation}\label{rem.jet.morphism2}
\phi_\square\colon J^r(E,n) \rTo \eta_r\,,\quad [L]^r_p  \mto
\left([L]^r_p\,,\,\square(L)(p) \right)
\end{equation}
The correspondence $\phi\mto \square_\phi$ is bijective.

\smallskip
If $\varrho\colon F\rTo M$ is a bundle, the \emph{$r$-jet bundle}
$J^r\varrho$ is defined as the $r$-jet of submanifolds which are
the image of local sections of $\varrho$. In the rest of the paper
$M$ will be an $n$-dimensional manifold. If $(x^\lambda)$ is a
chart on $M$, and $u^A$ fibred coordinates on $F$, then we denote
by $(x^\lambda,u^A_{x^{\bsi}})$ a chart on $J^r\varrho$, where the
subscript $x^{\bsi}$ is put to distinguish derivative coordinates
of a jet bundle from those of a jet of submanifolds. Roughly
speaking, we can interpret the manifold $M$ as the \emph{space of
parameters}. We can reproduce all constructions and definitions
related to jet of submanifolds in the case of jet bundle, then we
omit them.

\smallskip
Now we recall \cite{MaThesis,MaAnkara,MaVi02,MoVi94} the
construction of certain bundles which play a key role in the rest
of the paper. Let $r\geq 0$.

Let us consider the following bundles over $J^{r+1}(E,n)$: the
pull-back bundle
\[
T^{r+1,r}=\pi^*_{r+1,r}(TJ^r(E,n))
\]
of $TJ^r(E,n)$ through the map $\pi_{r+1,r}$, the sub-bundle
\[
H^{r+1,r}=\{(\theta_{r+1},v)\,\,/\,\,\theta_{r+1}\in J^{r+1}(E,n)
\,,\,v\in R_{\theta_{r+1}}\}
\]
of $T^{r+1,r}$, and the quotient bundle
\[
V^{r+1,r}=T^{r+1,r}\big{/}H^{r+1,r}.
\]

Of course, we have the following exact sequence:
\begin{equation}  \label{eq:contact}
\begin{diagram}
0 & \rTo & H^{r+1,r} & \rhTo^{D^{r+1}} & T^{r+1,r} & \rTo &
V^{r+1,r} & \rTo & 0\,,
\end{diagram}
\end{equation}
where $D^{r+1}$ is the natural inclusion.

\begin{remark}\label{rem.fibred.inclusion}
We can interpret $D^{r+1}$ as an inclusion
$J^{r+1}(E,n)\hookrightarrow
H^{r+1,r*}\otimes_{J^{r+1}(E,n)}T^{r+1,r}$.
\end{remark}

\smallskip
Also, let us introduce the exact sequence of bundles:
\begin{equation}
\begin{diagram}
0 & \rTo & H^{r+1,0} & \rhTo & T^{r+1,0} & \rTo^{\omega^{r+1}} &
V^{r+1,0} & \rTo & 0
\end{diagram}
\end{equation}
where $H^{r+1,0}=\pi^*_{r+1,1}(H^{1,0})$,
$T^{r+1,0}=\pi^*_{r+1,1}(T^{1,0})$ and
$V^{r+1,0}=\pi^*_{r+1,1}(V^{1,0})$ and $\omega^{r+1}$ is the
natural quotient projection.

\smallskip
Clearly, $H^{r+1,0}\simeq H^{r+1,r}$.

\smallskip
Local bases of $C^\infty(J^{r+1}(E,n))$-modules of local sections
of $H^{r+1,r}$, $H^{r+1,r\,*}$, $V^{r+1,0}$ and $V^{r+1,0\,*}$ are
respectively:
\begin{eqnarray*}
  D^{r+1,r}_\lambda&=&\frac{\partial }{\partial u^{\lambda }}+
  u_{\bsi,\lambda}^{j}\frac{\partial }{\partial u_{\bsi}^{j}}\,,
  \,\,\, \left[du^\lambda\right]=D^{r+1,r*}\circ(\pi^*_{r+1,r}(du^\lambda)),
  \\
  \left[\partial_{u^j}\right]&=&\omega^{r+1}\circ\partial_{u^j}
  \circ \pi_{r+1,0}\,, \,\,\, \omega^j=du^j-u^j_\lambda
  du^\lambda\,,
\end{eqnarray*}
where $|\bsi|\leq r$, and the pair $\bsi,\lambda$ stands for
$(\sigma_1,\dots,\sigma_k,\lambda)$.

%\\
%  \,,\quad \omega^j_{\bsi}=
%du_{\bsi}^{j}-u_{\bsi,\lambda}^{j}du^{\lambda},

\medskip
In the rest of the paper we denote
$\left[du^\lambda\right]$ and $\left[\partial_{u^j}\right]$
respectively by $du^\lambda$ and $\partial_{u^j}$.
\begin{remark}\label{rem.pseudo.nel.caso.fibrato}
In the case that $E$ is endowed with a bundle structure
$\varrho\colon E\rTo M$, we have that $H^{r+1,r}\simeq
\varrho^*_{r+1}(TM)$ and $V^{r+1,r}\simeq \pi^*_{r+1,r}(\ker
T\varrho_{r})$, where
$\varrho_k=\varrho\circ\varrho_{1,0}\circ\dots\circ\varrho_{k,k-1}$.
\end{remark}

\subsection{Symmetries of Differential Equations}

A \emph{classical external symmetry} of the equation
$\mathcal{E}\subset J^r(E,n)$ is a contact transformation of
$J^r(E,n)$ which preserves $\cE$. A \emph{classical internal
symmetry} of the equation $\mathcal{E}\subset J^r(E,n)$ is a
diffeomorphism of $\cE$ which preserves the induced Cartan
distribution $\cC(\cE)=\cC^r\cap T\cE$ on $\cE$. Analogously we
have an infinitesimal version of the previous definitions. We
would like to stress the we can speak about classical symmetries
even if that equation has no solutions. Classical symmetries send
solutions into solutions when these exist.

\smallskip
The \emph{1-prolongation} $\cE^1$ of an equation $\cE\subset
J^r(E,n)$ is the set
\[
\cE^1=\{\theta_{r+1}\in J^{r+1}(E,n)\,|\, \theta_r\in\cE,\,\,
R_{\theta_{r+1}}\subset T_{\theta_r}\cE\}.
\]
By iteration we can define the $l$-prolongation $\cE^l$. We denote
by $\theta_r^{(l)}$ the $l$-prolongation of the point
$\theta_r\in\cE$. The equation $\cE$ is said to be \emph{formally
integrable} if all the prolongations $\cE^l$ are smooth manifolds
and the maps $\pi_{r+l+1,r+l}|_{\cE^{l+1}}$ are fibre bundles.

In the rest of the section we deal just with equations which are
formally integrable.

A \emph{covering} is a pair $(\tilde{\cE},\tau)$, where
$\tilde{\cE}$ is a manifold provided with an $n$-dimensional
integrable distribution $\tilde{\cC}$\, and
$\tau\colon\tilde{\cE}\rTo\cE^\infty$ is a surjection such that
for any point $\tilde{\theta}\in\tilde{\cE}$ the tangent map
$T_{\tilde{\theta}}\tau$ maps isomorphically
$\cC_{\tilde{\theta}}(\tilde{\cE}^\infty)$ to the Cartan plane
$\cC_{\tau(\tilde{\theta})}(\cE^\infty)$. If
$\tilde{\cE}=\cE^\prime{}^\infty$ for some equation $\cE^\prime$
then the covering is called a \emph{covering equation}.

A vector field which belongs to $\cC^{\infty}(D)$ is called
\emph{trivial} as it is tangent to all the integral manifolds of
$\cC^{\infty}$. For this reason we call the quotient algebra
$\Sym\byd D_{\cC^{\infty}}/\cC^{\infty}(D)$ the algebra of
\emph{non-trivial symmetries} of the distribution $\cC^{\infty}$.

An element $X$ of $D_{\cC^{\infty}}$ is an \emph{infinitesimal
higher external symmetry} of the equation $\cE$ if $X$ is tangent
to $\cE^\infty$, and the equivalence class $[X]$ of $\Sym$ is
called a \emph{non-trivial infinitesimal higher external symmetry}
of $\cE$. The set of such classes is denoted by
$\Sym_{\text{ext}}(\cE)$.

Now, let us restrict our attention to $\cE^\infty$. An
\emph{infinitesimal higher internal symmetry} of the equation
$\cE$ is a symmetry of $\cC(\cE^{\infty})$. The set of such
symmetries is denoted by $D_{\cC}(\cE^{\infty})$. We denote by
$\cC D(\cE^{\infty})$ the ideal of vector fields on $\cE^\infty$
tangent to $\cC(\cE^{\infty})$, and we call \emph{trivial
symmetries} such fields. We denote by $\Sym_{\inter}(\cE)$ the
algebra $D_{\cC}(\cE^{\infty})\big/\cC D(\cE^{\infty})$ of
non-trivial infinitesimal higher internal symmetries of the
equation $\cE$.

The restriction map $\Sym_{\ext}(\cE)\rTo \Sym_{\inter}(\cE)$ is
surjective. Despite this fact, the set of classical external
symmetries does not project on the set of classical internal
symmetries (see \cite{Many99} for a general reference and theorem
\ref{th.ext.symm.are.subgroup.of.internal} for a concrete
example).

\section{Grassmannian Equivalent Connections}\label{sec.geodesic.equation}

In this section we construct a connection $\dot{\Gamma}$ on
$\pi_{1,0}$ starting from a torsion-free linear connection
$\Gamma$ on the tangent bundle $TE\rTo E$. By using $\dot{\Gamma}$
we find a section $\ddot{\Gamma}\colon J^1(E,n)\rTo J^2(E,n)$
whose image coincides with the equation of totally geodesic
submanifolds. We call such a map the $n$-Grassmannian structure
associated with $\Gamma$, as it turns out to be a distinguished
representative of the equivalence class of connections admitting
the same equation of unparametrized totally geodesic submanifolds.
Then two connections will be $n$-Grassmannian equivalent if and
only if they induce the same $n$-Grassmannian structure. As
by-product, we recover in the case $n=1$ the geometry of
projectively equivalent connections and projective structures.
Finally we introduce the equation of parametrized totally geodesic
submanifolds as a submanifold of $J^2{\pro_{M}}$ where
${\pro_{M}}\colon M\times E\rTo M$. Analogously to the
unparametrized case, we construct a connection
$\dot{\Gamma}_{\pro_{M}}$ on ${\pro_{M}}_{1,0}$ and an operator
$\ddot{\Gamma}_{\pro_{M}}$. We shall investigate in section
\ref{sec.covering} the relations between all the structures
constructed here.

Below, we fix some notations.

\smallskip
Let $\rho\colon F\rTo N$ be a bundle. A connection $\Gamma$ on
$\rho$ is a section of the bundle $\rho_{1,0}$. It can be
equivalently seen as tangent valued form and as a vertical valued
form, respectively
\begin{equation*}%\label{eq.tangent.and.vertical.valued.forms}
\Gamma\colon F\rTo \rho^*(T^*N)\otimes_{F}TF\,\,\,\,,\,\,\,
\upsilon_{\Gamma}\colon F\rTo T^*F\otimes_{F}VF\,,
\end{equation*}
where $\upsilon_\Gamma=\id_{TF}-\Gamma\circ T\rho$.

Let $(u^1,u^2\dots,u^{s})$ be a chart on $N$ and
$(u^1,u^2\dots,u^{s},\eta^1,\dots,\eta^l)$ a chart on $F$. We have
the local expressions
\begin{equation}\label{eq.local.expr.tangent.and.vertical.valued.forms}
\Gamma=du^{i}\otimes(\partial_{u^{i}}-\Gamma_i{}^k \,
\partial_{\eta^k}) \,\,\,\,,\,\,\, \upsilon_{\Gamma}=(d\eta^k +
\Gamma_i{}^k \, du^{i})\otimes\partial_{\eta^k} \,\,\,\,,\,\,\,
\Gamma_i{}^k\in C^{\infty}(F).
\end{equation}

If $\rho$ is a vector bundle, we can characterize the vertical
projection as a map $\upsilon_{\Gamma}\colon TF\rTo F$
\cite{Kling82,KMS93,Sau89}. Locally
\[
\upsilon_{\Gamma}= \left(d\eta^k+\Gamma_i{}^k{}_j \, \eta^{j} \,
du^{i}\right) \otimes p_k \,,\quad \Gamma_i{}^k{}_j \in
C^{\infty}(N)
\]
where $\{p_k\}$ is a local basis of $F$.

 The covariant derivative
$\nabla[\Gamma]_X$ of $\Gamma$ with respect to a vector field $X$
on $N$ acts on a section $s$ of $F$ by
$\nabla[\Gamma]_X(s)=\upsilon_\Gamma\circ Ts\circ X$.

For any map $f\colon P\rTo N$, we denote by $f^*(\Gamma)$ the
pull-back connection on the pull-back bundle $f^*(\rho)$.

\subsection{Totally Geodesic Submanifolds in the Jet of
Submanifolds}\label{sec.Totally.Geod.Subm.in.the.Jet.of.Subm}

From now on, we mean by $\Gamma$ a torsion-free linear connection
on the tangent bundle $TE\rTo E$. We denote by $\odot$ the
symmetric product. Let us introduce the operator
\begin{equation}\label{eq.operator.of.tot.geod.sub}
\operatorname{II}\colon J^2(E,n) \rTo H^{2,0\,\ast}\underset{J^{2}(E,n)}{\odot }%
H^{2,0\,\ast}\underset{J^{2}(E,n)}{\otimes }V^{2,0}
\end{equation}
defined by:
\begin{equation}\label{eq.operator.of.tot.geod.sub.2}
\operatorname{II}\,(X,Y,\omega)=\omega
\left(\nabla\left[\pi_{2,0}^*(\Gamma)\right]_{\widetilde{X}}(Y)\right)
\end{equation}
where $\widetilde{X}$ is a Cartan field on $J^2(E,n)$ projecting
on $X$ (see \eqref{eq.splitting.Cartan}). The previous definition
is well posed as the vertical part of $\tilde{X}$ gives no
contribution. We call $\operatorname*{II}$ the \emph{universal
second fundamental form} associated with $\Gamma$. The name is
justified by the fact that the restriction of $\operatorname*{II}$
on a $n$-dimensional submanifold $L$ of $E$ is the second
fundamental form on $L$. Then, taking into consideration
\eqref{rem.jet.morphism1} and \eqref{rem.jet.morphism2}, we give
the following
\begin{definition}
The equation $\cS$ of unparametrized totally geodesic submanifolds
is the submanifold of $J^2(E,n)$ given by $\operatorname*{II}=0$.
\end{definition}
Locally, such a submanifold is described by
\begin{multline}\label{eq.tot.geod.unparam}
u^k_{\lambda\xi}+\Gamma_\lambda{}^k{}_\xi + \Gamma_\lambda{}^k{}_i
\, u^i_\xi + \Gamma_j{}^k{}_\xi \, u^j_\lambda + \Gamma_j{}^k{}_i
\, u^j_\lambda u^i_\xi\\
-u^k_\beta\left(\Gamma_\lambda{}^\beta{}_\xi +
\Gamma_\lambda{}^\beta{}_i \, u^i_\xi + \Gamma_j{}^\beta{}_\xi \,
u^j_\lambda + \Gamma_j{}^\beta{}_i \, u^j_\lambda
u^i_\xi\right)=0.
\end{multline}

We could define the operator \eqref{eq.operator.of.tot.geod.sub}
by using the infinite jet of submanifolds $J^\infty(E,n)$. More
precisely we could consider the connection
$\pi_{\infty,0}^*(\Gamma)$, then define
\eqref{eq.operator.of.tot.geod.sub.2} simply by restriction on the
pseudo-horizontal bundle, as the pseudo-horizontal bundle
associated with $J^\infty(E,n)$ is a sub-bundle of
$TJ^\infty(E,n)$ and realize that the operator that we obtain is a
second order differential operator by a direct calculation.

Let us note that equation \eqref{eq.tot.geod.unparam} coincides,
in the case $n=1$, with the equation of unparametrized geodesics
\cite{Aminova01,MaThesis,MaAnkara}.

Also, \eqref{eq.tot.geod.unparam} can be solved with respect to
second derivatives. This suggests to think of this equation as the
image of some section of $\pi_{2,1}$. Below we propose an
intrinsic way to construct such a section.

\begin{lemma}\label{lemmma.isom.verticale}
$\ker T\pi_{1,0}\simeq H^{1,0\,*}\otimes_{J^1(E,n)}V^{1,0}$.
\end{lemma}
\begin{proof}
Let $v_{\theta_1}\in \ker T_{\theta_1}\pi_{1,0}$. Let us consider
a curve $\gamma (t)$ in the fibre
$\pi_{1,0}^{-1}(\pi_{1,0}(\theta_1))$ such that $\gamma
(0)=\theta_1$ and $\dot{\gamma}(0)=v_{\theta_1}$. Taking into
account remark \ref{rem.fibred.inclusion}, the above isomorphism
is realized by the map $v_{\theta_1} \mto \left.
\frac{d}{dt}\right| _{0}D^{1}(\gamma (t))$. In coordinates
\begin{equation*}
D^1(\gamma(t))=du^\lambda\otimes\left(\partial_{u^\lambda}%
+u^j_{\lambda}(t)\partial_{u^j}\right)\,,
\end{equation*}
then
\begin{equation*}
\left. \frac{d}{dt}\right| _{0}D^{1}(\gamma(t))=
\dot{u}^j_{\lambda}(0) \left(du^\lambda\otimes \partial_{u^j}%
\right).
\end{equation*}
\end{proof}

\begin{theorem}\label{theo.principale}
With the connection $\Gamma$ it is naturally associated a
connection $\dot{\Gamma}$ on $\pi_{1,0}$.
\end{theorem}
\begin{proof}
By previous lemma, the connection $\dot{\Gamma}$ can be
characterized as the projection
\begin{equation*}%\label{eq.isom.verticale}
\upsilon_{\dot{\Gamma}}\colon TJ^1(E,n)\rTo
H^{1,0\,*}\underset{J^1(E,n)}\otimes V^{1,0}.
\end{equation*}
This projection is defined by the following commutative diagram:
\begin{diagram}
TJ^1(E,n) &
\rTo^{\upsilon_{\dot{\Gamma}}}  & H^{1,0\,*}\underset{J^1(E,n)}\otimes V^{1,0} \\
\dTo^{TD^1}  &  & \uTo_{\id_{H^{1,0\,*}}\otimes\,\omega^{1}} \\
T\left(H^{1,0\,*}\underset{J^1(E,n)}\otimes T^{1,0}\right) &
\rTo^{\upsilon_{\Omega}} & {H^{1,0\,*}\underset{J^1(E,n)}\otimes
T^{1,0}}
\end{diagram}
where $\Omega=\Xi\otimes\pi^*_{1,0}(\Gamma)$ and $\Xi$ is a linear
connection on the bundle $H^{1,0\,*}$. Let
$(u^\lambda,u^j,\eta_\lambda)$ be a chart on $H^{1,0\,*}$
associated with $du^\lambda$. Then:
\[
\upsilon_{\Xi}=\left(d\eta_\xi + \Xi_{A}{}_{\xi}{}^{\lambda} \,
\eta_\lambda \, du^A + \Xi^\alpha_h{\,}_{\xi}{}^{\lambda} \,
\eta_\lambda \, du^h_\alpha\right)\otimes du^\xi.
\]

We continue to denote by $\Gamma$ the connection
$\pi^*_{1,0}(\Gamma)$. We shall see that the map on the top row of
the above diagram gives a connection on $J^1(E,n)\rTo E$ which is
independent of the connection $\Xi$.

Let $(u^\lambda,u^j,u^j_\lambda,\eta_{\lambda}^A)$ be the local
chart on $H^{1,0\,*}\otimes_{J^1(E,n)} T^{1,0}$ associated with
the basis $\{du^\lambda\otimes \partial_{u^A}\}$. We have that

\begin{equation*}
TD^{1}=du^A\otimes\partial_{u^A}+du^j_\lambda\otimes
\left(\partial_{u^j_\lambda}+\partial_{\eta_\lambda^j}\right)
\end{equation*}
and
\[
\upsilon_{\Omega} = \left( d\eta^C_\xi +
\Omega_A{\,}^C_\xi{\,}^\lambda_B \, \eta^B_\lambda +
\Omega^\alpha_h{\,}^C_\xi{\,}^\lambda_B \, du^h_\alpha \right)
\otimes du^\xi\otimes\partial_{u^C}
\]
with
\begin{equation}\label{eq.simboli.Christo.Omega}
\Omega_A{\,}^C_\xi{\,}^\lambda_B =
\delta_B^C\,\Xi_{A}{}_{\xi}{}^{\lambda}+
\delta_\xi^\lambda\Gamma_A{}^C{}_B\,, \quad
\Omega^\alpha_h{\,}^C_\xi{\,}^\lambda_B =
\delta_B^C\,\Xi^\alpha_h{}_\xi{}^{\lambda}+
\delta_\xi^\lambda\Gamma^\alpha_h{}\,^C{}_{B}=\delta_B^C\,\Xi^\alpha_h{}_\xi{}^{\lambda}.
\end{equation}
Then
\begin{eqnarray*}
\upsilon _{\Omega }\circ TD^{1} &=& \big(
\Omega_A{\,}^C_\xi{\,}^\lambda_\lambda \, du^{A} +
\Omega_A{\,}^C_\xi{\,}^\lambda_i \, u^i_\lambda \,du^{A} + \\
&&\,+\,\delta^{C}_{h}\,\delta^{\alpha}_{\xi }\,du_{\alpha
}^{h}+\Omega^\alpha_h{\,}^C_\xi{\,}^\lambda_\lambda \, du_{\alpha
}^{h} + \Omega^\alpha_h{\,}^C_\xi{\,}^\lambda_i \, u_{\lambda
}^{i} \, du_{\alpha }^{h} \big)\otimes du^{\xi }\otimes \partial
_{u^{C}}.
\end{eqnarray*}

Now it remains to project the previous expression on
$H^{1,0\,*}\otimes V^{1,0}$. We obtain
\begin{eqnarray*}
\upsilon_{\dot{\Gamma}} &=& \left( \left(
\Omega_A{\,}^k_\xi{\,}^\lambda_\lambda +
\Omega_A{\,}^k_\xi{\,}^\lambda_i \, u^i_\lambda - u^k_\beta \left(
\Omega_A{\,}^\beta_\xi{\,}^\lambda_\lambda \,+
\Omega_A{\,}^\beta_\xi{\,}^\lambda_i \, u^i_\lambda \right)
\right)du^A + \right.
\\
&+& \left. \left(\delta_h^k\delta_\xi^\alpha +
\Omega^\alpha_h{\,}^k_\xi{\,}^\lambda_\lambda +
\Omega^\alpha_h{\,}^k_\xi{\,}^\lambda_i \, u^i_\lambda -
u^k_\beta\left(\delta_h^\beta\delta_\xi^\alpha +
\Omega^\alpha_h{\,}^\beta_\xi{\,}^\lambda_\lambda +
\Omega^\alpha_h{\,}^\beta_\xi{\,}^\lambda_i \, u^i_\lambda\right)
\right)du^h_\alpha \right)\otimes du^\xi\otimes\partial_{u^k}.
\end{eqnarray*}
In view of \eqref{eq.simboli.Christo.Omega},
$\upsilon_{\dot{\Gamma}}$ does not depend on $\Xi$. In fact, on
one hand we have
\[
\Omega_A{\,}^k_\xi{\,}^\lambda_\lambda +
\Omega_A{\,}^k_\xi{\,}^\lambda_i \, u^i_\lambda - u^k_\beta \left(
\Omega_A{\,}^\beta_\xi{\,}^\lambda_\lambda \,+
\Omega_A{\,}^\beta_\xi{\,}^\lambda_i \, u^i_\lambda \right) =
\Gamma_{A}{}^{k}{}_{\xi} + \Gamma_{A}{}^{k}{}_{i} \, u^i_\xi-
u^k_\beta ( \Gamma_{A}{}^{\beta}{}_{\xi} +
\Gamma_{A}{}^{\beta}{}_{i} \, u^i_\xi ).
\]
On the other hand we have
\[
\delta_h^k\delta_\xi^\alpha +
\Omega^\alpha_h{\,}^k_\xi{\,}^\lambda_\lambda +
\Omega^\alpha_h{\,}^k_\xi{\,}^\lambda_i \, u^i_\lambda -
u^k_\beta\left(\delta_h^\beta\delta_\xi^\alpha +
\Omega^\alpha_h{\,}^\beta_\xi{\,}^\lambda_\lambda +
\Omega^\alpha_h{\,}^\beta_\xi{\,}^\lambda_i \,
u^i_\lambda\right)=\delta_h^k\delta_\xi^\alpha.
\]

Taking into consideration that the connection $\dot{\Gamma}$ is
characterized by
\begin{equation}\label{eq.local.expr.of.upsilon.dot.nabla}
\upsilon_{\dot{\Gamma}}=\left(du^k_\xi+\dot{\Gamma}_A{}^k{}_\xi \,
du^A\right) \otimes du^\xi\otimes\partial_{u^k}\,,
\end{equation}
the theorem follows by putting
\[
\dot{\Gamma}_A{}^k{}_\xi=\Gamma_A{}^k{}_\xi + \Gamma_A{}^k{}_i \,
u^i_\xi- u^k_\beta(\Gamma_A{}^\beta{}_\xi + \Gamma_A{}^\beta{}_i
\, u^i_\xi).
\]
\end{proof}

\begin{corollary}
The correspondence $\Gamma\mto\dot{\Gamma}$ is injective.
\end{corollary}
\begin{proof}
It is a matter of computation, taking into consideration that
$\dot{\Gamma}_A{}^k{}_\xi$ are polynomials of second order degree
in the first derivatives.
\end{proof}

\begin{theorem}
With any connection $\dot{\Gamma}$ on $\pi_{1,0}$ it is naturally
associated a section $\ddot{\Gamma}$ of $\pi_{2,1}$.
\end{theorem}
\begin{proof}
Taking into account remark \ref{rem.fibred.inclusion}, we can
characterize a section $\ddot{\Gamma}$ of $\pi_{2,1}$ as the
unique section of $\pi_{2,1}$ such that
$D^2\circ\ddot{\Gamma}=D^1\lrcorner\,\dot{\Gamma}$.
\end{proof}

%
%In view of remark \ref{rem.bundle.on.r.rather.than.r+1}, we can
%interpret $D^2$ as the fibred inclusion $J^{2}(E,n)\hookrightarrow
% (H^{1,0})^*\otimes_{J^{1}(E,n)}TJ^{1}(E,n)$. Then
%we obtain a section of $\pi_{2,1}$ if we consider a section of
%$(H^{1,0}){}^*\otimes_{J^1(E,n)} TJ^1(E,n)$ and restrict it to the
%image of $D^2$. In order to obtain such a section it is sufficient
%to consider the map
%\[
%\ddot{\Gamma}=D^1\lrcorner\,\dot{\Gamma}.
%\]
%
The local expression of $D^2\circ\ddot{\Gamma}$ is given by
\begin{equation}\label{eq.local.exp.of.ddotnabla}
D^2\circ\ddot{\Gamma} =
du^\lambda\otimes\left(\partial_{u^\lambda} +
u^j_\lambda\partial_{j} - \left( \dot{\Gamma}_\lambda{}^k{}_\xi +
u^j_\lambda\,\dot{\Gamma}_j{}^k{}_\xi\right)\partial_{u^k_\xi}\right).
\end{equation}

\begin{proposition}
The image of $\,\ddot{\Gamma}$ coincides with the equation of
totally geodesic submanifolds \eqref{eq.tot.geod.unparam}.
\end{proposition}
\begin{proof}
It is sufficient to take into consideration expressions
\eqref{eq.tot.geod.unparam} and \eqref{eq.local.exp.of.ddotnabla}.
\end{proof}

\begin{remark}\label{rem.distribution}
Of course, we can interpret $D^2\circ\ddot{\Gamma}$ as a
distinguished distribution on $J^1(E,n)$. More precisely
$D^2(\ddot{\Gamma}(\theta_1))\equiv R_{\ddot{\Gamma}(\theta_1)}$,
where $R_{\ddot{\Gamma}(\theta_1)}$ is the $R$-plane associated
with the point $\ddot{\Gamma}(\theta_1)$. From now on we shall
denote by $R\circ\ddot{\Gamma}$ such a distribution.
\end{remark}

\begin{remark}
The papers \cite{JanyModyGalilei,JanyModyRoma} are devoted to the
formulation of Galilean relativistic mechanics based on jet spaces
and cosymplectic forms. Further developments of this theory are
present \cite{JanyModyGeneral}, where the theory is extended to
Einstein relativistic mechanics. In \cite{MaVi04} variational
aspects are studied. Theorem \ref{theo.principale} is inspired by
\cite{JanyModyGeneral}, but the results obtained here are much
more general. In fact in \cite{JanyModyGeneral} $\dot{\Gamma}$ and
$\ddot{\Gamma}$ have been constructed in the case of jets of
time-like curves of an oriented manifold. Because of orientation,
the pseudo-horizontal bundle turns out to be trivial, and a
connection on it with vanishing Christoffel symbols always exists.
Here we construct the connection $\dot{\Gamma}$ in the case of
$n$-dimensional submanifolds without requiring any orientation,
and we prove that $\dot{\Gamma}$ is independent of the connection
$\Xi$ on $H^{1,0\,*}$. This gives a natural character to
$\dot{\Gamma}$. Also, from a physical viewpoint, taking into
consideration $\ddot{\Gamma}$, we obtain as by-product, in the
case $n=1$ and for time-like curves, the equation of motion of one
relativistic particle \cite{JanyModyGeneral,MaVi04}.
\end{remark}

\begin{definition}
We say that two torsion-free linear connections on $TE\rTo E$ are
\emph{$n$-Grassmannian equivalent} if they have the same equation
of $n$-dimensional totally geodesic submanifolds.
\end{definition}

We note that in the case $n=1$ the distribution
$R\circ\ddot{\Gamma}$ is integrable and the above definition
coincides with that of projectively equivalent connections.

\smallskip
Our reasoning leads naturally to the following
\begin{theorem}\label{th.equiv.grass.conn}
Two torsion-free linear connections $\Gamma$ and ${\Gamma^\prime}$
are $n$-Grassmannian equivalent if and only if
$\,\ddot{\Gamma}=\ddot{{\Gamma^\prime}}$.
\end{theorem}

\begin{proposition}
The following quantities
\begin{equation}\label{eq.Grassmannian.invariants}
\Gamma_\lambda{}^k{}_\xi\,,
%%%
\quad  \delta_\xi^\beta \, \Gamma_\lambda{}^k{}_i +
\delta_\lambda^\beta \, \Gamma_i{}^k{}_\xi - \delta_i^k \,
\Gamma_\lambda{}^\beta{}_\xi\,,
%%%
\quad \delta_\lambda^\alpha \, \delta_\xi^\beta \,
\Gamma_j{}^k{}_i - \delta^\alpha_\xi \, \delta^k_i \,
\Gamma_\lambda{}^\beta{}_j - \delta_\lambda^\alpha \, \delta_i^k
\, \Gamma_j{}^\beta{}_\xi\,,
%%%
\quad \Gamma_j{}^\beta{}_i
\end{equation}
are $n$-Grassmannian invariants, that is they do not change for a
connection which is $n$-Grassmannian equivalent to $\Gamma$.
\end{proposition}
\begin{proof}
Equation \eqref{eq.tot.geod.unparam} is a system of differential
equations $u^j_{\lambda\xi}+P^j_{\lambda\xi}=0$ where
$P^j_{\lambda\xi}$ are polynomials of third order in the first
derivatives. Quantities \eqref{eq.Grassmannian.invariants} are
respectively the coefficients of order $0$, $1$, $2$, and $3$ of
such a system.
\end{proof}

\begin{corollary}
Two torsion-free linear connections $\Gamma$ and ${\Gamma^\prime}$
are $n$-Grassmannian equivalent if they have the same
$n$-Grassmannian invariants. Locally they are related by:
%\begin{equation}\label{eq.grass.equiv.conn}
\begin{eqnarray}\label{eq.grass.equiv.conn}
\Gamma_\lambda{}^k{}_\xi & = & {\Gamma'}_\lambda{}^k{}_\xi \nonumber \\
%%\\
\delta^\beta_\xi\left(\Gamma_\lambda{}^k{}_i-{\Gamma'}_\lambda{}^k{}_i\right)
+
\delta^\beta_\lambda\left(\Gamma_i{}^k{}_\xi-{\Gamma'}_i{}^k{}_\xi\right)
& = &
\delta_i^k\left(\Gamma_\lambda{}^\beta{}_\xi-{\Gamma'}_\lambda{}^\beta{}_\xi\right) \nonumber \\
%%\\
2\delta^\beta_\lambda\delta^\alpha_\xi\left(\Gamma_j{}^k{}_i-{\Gamma'}_j{}^k{}_i\right)
& = &
\delta^\alpha_\xi\delta^k_j\left(\Gamma_\lambda{}^\beta{}_i-{\Gamma'}_\lambda{}^\beta{}_i\right)
+ \delta^k_j\delta^\alpha_\lambda\left(\Gamma_i{}^\beta{}_\xi-{\Gamma'}_i{}^\beta{}_\xi\right)+\\
&  &
\delta^\alpha_\xi\delta_i^k\left(\Gamma_\lambda{}^\beta{}_j-{\Gamma'}_\lambda{}^\beta{}_j\right)
+\delta^k_i\delta^\alpha_\lambda\left(\Gamma_j{}^\beta{}_\xi-{\Gamma'}_j{}^\beta{}_\xi\right)
\nonumber
\\
%%\\
\Gamma_j{}^\beta{}_i & = & {\Gamma'}_j{}^\beta{}_i \nonumber
\end{eqnarray}
%\end{equation}
\end{corollary}
\begin{proof}
Straightforward.
\end{proof}

\begin{definition}
We call $\ddot{\Gamma}$ the \emph{$n$-Grassmannian structure}
associated with the connection $\Gamma$.
\end{definition}

For $n=1$ the definition of $1$-Grassmannian structure can be
interpreted as the projective structure of $E$ associated with the
connection $\Gamma$. More precisely we have the following

\begin{proposition}[\cite{Aminova01}]
Let $n=1$. Then quantities \eqref{eq.Grassmannian.invariants} can
be expressed in terms of projective invariants
\eqref{eq.projective.invariants}.
\end{proposition}

\begin{remark}\label{re.proj.structure.and.proj.eq.conn}
In \cite{Aminova01} the author obtains the result of the above
proposition by considerations of local character. Here such result
emerges naturally from our geometrical setting. Also, the
advantage of considering the projective structure as a particular
distribution on $J^1(E,1)$ (see also remark
\ref{rem.distribution}) rather than a particular subbundle of the
second frame bundle of $E$ is that of giving, for instance, a
clear interpretation of the relationship between projective and
contact projective symmetries (see section
\ref{sec.symm.of.geodesic.equation}).
\end{remark}

As a natural consequence of our reasoning we have the following
\begin{proposition}
Let $n=1$. Then equations \eqref{eq.grass.equiv.conn} reduce to
\eqref{eq.proj.equiv.connecions.locally}.
\end{proposition}

\subsection{Totally Geodesic Submanifolds in the Jet of a Trivial
Bundle}\label{sec.Totally.Geod.Subm.in.the.Jet.Bundle}

Let $M$ be an $n$-dimensional manifold. Let ${\pro_{M}}\colon
M\times E\rTo M$ be the projection on the first factor and
${\pro_{E}}\colon M\times E\rTo E$ the projection on the second
factor. We recall the projections ${\pro_M}_{r,r-1}\colon
J^{r}\pro_M\rTo J^{r-1}\pro_M$ and ${\pro_M}_r\colon
J^{r}\pro_M\rTo M$ where
${\pro_M}_r={\pro_M}\circ{\pro_M}_{1,0}\circ\dots\circ{\pro_M}_{r,r-1}$.
Let $\left(x^\lambda,u^A_{x^{\bsi}}\right)$ be a chart on
$J^r{\pro_{M}}$, $|\bsi|\leq r$. A parametrized $n$-dimensional
submanifold is an immersion $s\colon U\subset M\rTo E$. We notice
that we can interpret $s$ as a local section of ${\pro_{M}}$.
\begin{definition}
We denote by $\Imm J^r{\pro_{M}}$ the subset of $J^r{\pro_{M}}$ of
the $r$-jets of immersions of $M$ into $E$.
\end{definition}
We notice that $\Imm J^r{\pro_{M}}$ is an open dense submanifold
of $J^r{\pro_{M}}$.

Let $\Theta$ be a torsion-free linear connection on $TM \rTo M$.
Then $\Gamma\oplus\Theta$ is a connection on $TJ^0{\pro_{M}}\rTo
J^0{\pro_{M}}$.

Now, taking into consideration that $V^{2,0}\simeq
\pro_{M_{2,0}}^*(\pro_E^*(TE))$, and following the same reasoning
as section \ref{sec.Totally.Geod.Subm.in.the.Jet.of.Subm}, we can
define the following operator:
\[
\operatorname{II}_{\pro_{M}}\colon \Imm
J^2{\pro_{M}}\longrightarrow
 {\pro^*_{M}}_2(T^*M)\underset{J^2{\pro_{M}}}\odot {\pro^*_{M}}_2
(T^*M)\underset{J^2{\pro_{M}}}\otimes
{\pro^*_M}_{2,0}(\pro_E^*(TE)).
\]

\begin{definition}
The equation $\cS_{\pro_{M}}$ of parametrized totally geodesic
submanifolds is the submanifold of $\Imm J^2{\pro_{M}}$ given by
$\operatorname*{II}_{\pro_{M}} =0$.
\end{definition}

Locally, a parametrized totally geodesic submanifold is described
by
\begin{equation}\label{eq.tot.geod.param}
u^C_{x^\xi x^\lambda}+\Gamma_A{}^C{}_B \, u^A_{x^\xi}
u^B_{x^\lambda} - \Theta_\xi{}^\eta{}_\lambda \, u^C_{x^\eta} =0.
\end{equation}

\begin{proposition}
With the pair ($\Gamma,\Theta$) it is associated a connection
$\dot{\Gamma}_{\pro_{M}}$ on ${\pro_{M}}_{1,0}$ and a section
$\ddot{\Gamma}_{\pro_{M}}$ of ${\pro_{M}}_{2,1}$ whose image is
the equation of parametrized totally geodesic submanifolds.
\end{proposition}
\begin{proof}
We have that $J^1{\pro_{M}}\simeq
{\pro^*_{E}}(TE)\otimes{\pro^*_{M}}(T^*{M})$ \cite{KMS93}. Then we
can consider the connection
$\dot{\Gamma}_{\pro_{M}}={\pro^*_{E}}(\Gamma)\otimes{\pro^*_{M}}(\Theta^*)$,
where $\Theta^*$ is the dual connection of $\Theta$. The local
expression of $\dot{\Gamma}_{\pro_{M}}$ as tangent valued form is:
\begin{equation}\label{eq.local.expr.of.dot.nabla.varrho}
\dot{\Gamma}_{\pro_{M}}=
dx^\lambda\otimes\left(\partial_{x^\lambda} +
\Theta_\lambda{}^\xi{}_\eta \, u^A_{x^\xi}
\partial_{u^A_{x^\eta}} \right) +
du^B\otimes\left(\partial_{u^B}-\Gamma_B{}^C{}_A \,
u^A_{x^\xi}\partial_{u^C_{x^\xi}}\right)
\end{equation}

Similarly to the reasoning adopted in the previous section, we can
associate with $\dot{\Gamma}_{\pro_{M}}$ a section
$\ddot{\Gamma}_{\pro_{M}}$ of ${\pro_{M}}_{2,1}$. The image of
such a section is described by \eqref{eq.tot.geod.param}.
\end{proof}

\begin{remark}
If we consider the case $M=\R$, then equation
\eqref{eq.tot.geod.param} is the equation of (parametrized)
geodesics. Also, the connection $\Theta$ is flat, and then we can
find a system of coordinates where the Christoffel symbol
$\Theta_0{}^0{}_0$ vanishes. This change of coordinates is
equivalent to introduce an affine parameter.
\end{remark}

\section{Symmetries of the Equation of Totally Geodesic Submanifolds}\label{sec.symm.of.geodesic.equation}
In this section we discuss the symmetries of both the equations
$\cS$ and $\cS_{\pro_{M}}$, by generalizing results obtained in
\cite{MaAnkara}. We see, in this case, how classical internal
symmetries are the most general symmetries. We recall that with
$\ddot{\Gamma}$ we can associate the distribution
$R\circ\ddot{\Gamma}$ (see remark \ref{rem.distribution}). The
same consideration holds true for $\ddot{\Gamma}_{\pro_{M}}$. We
call $\cS$ (resp. $\cS_{\pro_{M}}$) \emph{integrable} if the
distribution $R\circ\ddot{\Gamma}$ (resp.
$R\circ{\ddot{\Gamma}_{\pro_{M}}}$) is integrable. We note that we
can speak about symmetries of $\cS$ (resp. $\cS_{\pro_{M}}$) even
if it has no solutions. They are symmetries of the
(non-integrable) distribution $R\circ\ddot{\Gamma}$ (resp.
$R\circ{\ddot{\Gamma}_{\pro_{M}}}$).

\begin{theorem}\label{th.formally.integrable}
Let $\cS$ (resp. $\cS_{\pro_{M}}$) be integrable. Then all the
$l$-prolongations $\cS^l$ (resp. $\cS^l_{\pro_{M}}$) of $\cS$
(resp. $\cS_{\pro_{M}}$) are diffeomorphic to $J^1(E,n)$ (resp.
$J^1{\pro_{M}}$), and the induced Cartan distribution $\cC(\cS^l)$
(resp. $\cC(\cS^l_{\pro_{M}})$) is isomorphic to
$R\circ{\ddot{\Gamma}}$ (resp.
$R\circ{\ddot{\Gamma}_{\pro_{M}}}$).
\end{theorem}
\begin{proof}
We give the proof for the equation $\cS$. The same reasoning holds
true for the equation $\cS_{\pro_{M}}$.

\smallskip
Let $\theta_1\in J^1(E,n)$. For each point
$\ddot{\Gamma}(\theta_1)$ there exists only one $R$-plane which
contains such a point and which is contained in
$T_{\ddot{\Gamma}(\theta_1)}\cS$. In fact if we represent the
point $\ddot{\Gamma}(\theta_1)$ as the pair
$(\theta_1,R_{\ddot{\Gamma}(\theta_1)})$, then such an $R$-plane
is given by
$T_{\theta_1}\ddot{\Gamma}(R_{\ddot{\Gamma}(\theta_1)})$, and the
existence follows. On the other hand if two $R$-planes $R$ and
$\tilde{R}$ contain the same point $\ddot{\Gamma}(\theta_1)$, then
the difference between a vector $v\in R$ and a vector
$\tilde{v}\in\tilde{R}$ belongs to the kernel of
$T_{\ddot{\Gamma}(\theta_1)}\pi_{2,1}$. Of course we have that
$T_{\ddot{\Gamma}(\theta_1)}\cS\cap\ker
T_{\ddot{\Gamma}(\theta_1)}\pi_{2,1}=0$, and the uniqueness is
proved.

From this discussion we have that
\[
\cC_{\ddot{\Gamma}(\theta_1)}(\cS)=
\cC^2_{\ddot{\Gamma}(\theta_1)}\cap
T_{\ddot{\Gamma}(\theta_1)}\cS=
T_{\theta_1}\ddot{\Gamma}(R_{\ddot{\Gamma}(\theta_1)})\simeq
R_{\ddot{\Gamma}(\theta_1)}.
\]

We notice that $(\ddot{\Gamma}(\theta_1))^{(1)}$ is the point
characterized by the pair
$(\ddot{\Gamma}(\theta_1),T_{\theta_1}\ddot{\Gamma}(R_{\ddot{\Gamma}(\theta_1)}))$.
Also, we have
\[
\cC_{(\ddot{\Gamma}(\theta_1))^{(1)}}(\cS^1)=
T_{\theta_1}\ddot{\Gamma}^{(1)}(R_{\ddot{\Gamma}(\theta_1)})
\simeq R_{\ddot{\Gamma}(\theta_1)},
\]
where $\ddot{\Gamma}^{(1)}\colon\theta_1\mto
(\ddot{\Gamma}(\theta_1))^{(1)}$.  By iterating the previous
construction, we obtain that the $l$-prolongation of the equation
$\cS$ is diffeomorphic to $\cS$, and so to $J^1(E,n)$, and that
every Cartan plane $\cC_{(\ddot{\Gamma}(\theta_1))^{(l)}}(\cS^l)$
is isomorphic to $R_{\ddot{\Gamma}(\theta_1)}$.
\end{proof}

The section $\ddot{\Gamma}$ splits canonically the Cartan
distribution $\cC^1$ of $J^1(E,n)$. Namely, for each $\theta_1\in
J^1(E,n)$ we have that $
\cC^1_{\theta_1}=R_{\ddot{\Gamma}(\theta_1)}\oplus\ker
T_{\theta_1}\pi_{1,0}$ (see also remark \ref{rem.distribution}).
The sequence
\[
0\rTo H\rhTo TJ^1(E,n) \rTo^{W} TJ^1(E,n)/H \rTo 0\,,
\]
where $H=\{(\theta_1,v)\in TJ^1(E,n)\,\, |\,\, v\in
R_{\ddot{\Gamma}(\theta_1)}\}$ and $W$ is the natural projection,
is exact. We note that $H$ is nothing but the pull-back bundle of
$H^{2,1}$ through $\ddot{\Gamma}$.

The previous sequence induces the following sequence of modules of
sections:
\begin{equation*}%\label{eq.sequenza.nontrivialization.map}
0\rTo \cH\rhTo \chi(J^1(E,n)) \rTo^{\cW} \cV \rTo 0.
\end{equation*}
Similarly, we have the following exact sequence:
\begin{equation*}%\label{eq.sequenza.nontrivialization.map.bundle}
0\rTo \cH_{\pro_{M}}\rhTo \chi(J^1{\pro_{M}})
\rTo^{\cW_{\pro_{M}}} \cV_{\pro_{M}} \rTo 0,
\end{equation*}
where $\cH_{\pro_{M}}$ is the module of sections of the bundle
$H_{\pro_{M}}=\{(\tilde{\theta}_1,v)\in TJ^1{\pro_{M}}\,\, |\,\,
v\in R_{\ddot{\Gamma}_{\pro_{M}}(\tilde{\theta_1})}\}$ and
$\cW_{\pro_{M}}$ is the natural projection.

\begin{corollary}\label{cor.caratterizzazione.simmetrie.geodetiche}
Let $\cS$ (resp. $\cS_{\pro_{M}}$) be integrable. Then the algebra
of trivial internal symmetries of $\cS$ (resp. $\cS_{\pro_{M}}$)
is isomorphic to $\cH$ (resp. $\cH_{\pro_{M}}$). The algebra of
higher internal symmetries of $\cS$ (resp. $\cS_{\pro_{M}}$) is
isomorphic to the algebra of classical internal symmetries of
$\cS$ (resp. $\cS_{\pro_{M}}$), and coincides with the algebra of
vector fields on $J^1(E,n)$ (resp. $J^1{\pro_{M}}$) preserving the
$n$-dimensional distribution $R\circ\ddot{\Gamma}$ (resp.
$R\circ\ddot{\Gamma}_{{\pro_{M}}}$). The algebra of non-trivial
higher internal symmetries of $\cS$ (resp. $\cS_{\pro_{M}}$) is
the projection through $\cW$ (resp. $\cW_{\pro_{M}}$) of the
algebra of internal symmetries.
\end{corollary}
\begin{proof}
From Theorem \ref{th.formally.integrable} it follows that
$(\cS^\infty,\cC(\cS^\infty))$ is isomorphic to
$(J^1(E,n),R\circ\ddot{\Gamma})$. The same holds true for
$\cS_{\pro_{M}}$ and the proposition is proved.
\end{proof}
We would like to point out that even if $\cS$ (resp.
$\cS_{\pro_{M}}$) is not integrable, the results contained in the
previous corollary still remain true. The only difference is that
we can not speak about trivial symmetries or higher symmetries in
the sense of section \ref{sec.preliminary.definitions}. But the
characterization of classical internal symmetries is the same.

\begin{definition}
The group of external point symmetries of $\cS$ is called the
group of \emph{Grassmannian transformations} of $E$. The group of
internal symmetries of $\cS$ is called the group of \emph{contact
Grassmannian transformations} of $E$.
\end{definition}
\begin{remark}
We would like to note that in the case that $\dim E=n+1$, in view
of Lie-B\"{a}cklund theorem \cite{Many99,KrVe98,Vin81}, there are
external symmetries of $\cS$ which come from a contact
transformation of $J^1(E,n)$ but not from a transformation of $E$.
\end{remark}
\begin{proposition}
Let $\cS$ be integrable. An $n$-Grassmannian transformation sends
$n$-dimensional totally geodesic submanifolds into $n$-dimensional
totally geodesic submanifolds.
\end{proposition}
\begin{proof}
It is sufficient to apply the definitions.
\end{proof}

The following theorem generalizes a result contained in
\cite{MaThesis} and inspired by \cite{Many99}.
\begin{theorem}\label{th.ext.symm.are.subgroup.of.internal}
The group of classical external symmetries of $\cS$ (resp.
$\cS_{\pro_{M}}$) is a subgroup of the internal ones.
\end{theorem}
\begin{proof}
We give the proof in the case of the equation $\cS$. The result in
the case of the equation $\cS_{\pro_{M}}$ is attained by using a
similar reasoning.

We recall that any transformation $G\colon J^2(E,n)\rTo J^2(E,n)$
is a classical external symmetry of $\cS$ if $G$ is a contact
transformation and $G(\cS)=\cS$. Any such transformation projects
on a Lie transformation $G_{(1)}$ of $J^1(E,n)$. Also, we have
that $(G_{(1)})^{(1)}=G$. This means that $G$ is a Lie
transformation if and only if $G_{(1)}$ is a Lie transformation.
Moreover, the condition that $G$ preserves $\cS$ is equivalent to
require that $G_{(1)}$ preserves $R\circ\ddot{\Gamma}$. Then we
can identify any external symmetry with a diffeomorphism on
$J^1(E,n)$ which preserves both the Cartan distribution on
$J^1(E,n)$ and the distribution $R\circ\ddot{\Gamma}$.
\end{proof}
\begin{remark}
In the case $n=1$ the distribution $R\circ{\ddot{\Gamma}}$ (resp.
$R\circ{\ddot{\Gamma}_{\pro_{M}}}$) is of course integrable. In
this case, the infinitesimal point external symmetries of $\cS$
coincide with the projective fields on $E$, which have been
largely studied in literature (see \cite{Aminova01,Aminova04} and
the references therein). The previous theorem says that projective
symmetries form a sub-class of contact projective symmetries.
Studying such symmetries is interesting, for instance, in general
relativistic mechanics. In fact in \cite{MaVi04} the scheme of one
relativistic particle is constructed in the framework of jets of
submanifolds.
\end{remark}

\section{Covering of the Equation of Totally Geodesic Submanifolds}\label{sec.covering}
\begin{definition}
The \emph{r-prolongation} $\pro_{E}^{(r)}\colon \Imm
J^r{\pro_{M}}\rTo J^r(E,n)$ of the trivial projection
${\pro_{E}}\equiv \pro_{E}^{(0)}\colon M\times E\rTo E$ is defined
by $\pro_{E}^{(r)}\left([s]^r_{x}\right)=[{\pro_{E}}(\im
s)]^r_{{\pro_{E}}(s(x))}$, $x\in M$.
\end{definition}
Roughly speaking the projection $\pro_{E}^{(r)}$ sends jets of a
parametrized submanifold in jets of the same submanifold regarded
as an unparametrized submanifold.
\begin{lemma}
\begin{eqnarray}\label{equ.prolongations.of.tau}
\pro_{E}^{(r)}(\tilde{\theta}_{r})  & \equiv
\pro_{E}^{(r)}(\tilde{\theta
}_{r-1},R_{\tilde{\theta}_{r}})=\left(
\pro_{E}^{(r-1)}(\tilde{\theta
}_{r-1}),T_{\tilde{\theta}_{r-1}}\pro_{E}^{(r-1)}(R_{\tilde{\theta}_{r}%
})\right)  \\
& \equiv\left(  \pro_{E}^{(r-1)}(\tilde{\theta}_{r-1}),R_{\pro_{E}%
^{(r)}(\tilde{\theta}_{r})}\right) \nonumber  .
\end{eqnarray}
\end{lemma}
\begin{proof}
It is a straightforward computation.
\end{proof}
We have the following commutative diagram:
\begin{equation}\label{equ.comm.diagram.tau}
\begin{diagram}
\ldots & \rTo^{{\pro_{M}}_{r+1,r}} & \Imm J^r{\pro_{M}} &
\rTo^{{\pro_{M}}_{r,r-1}} & \ldots & \rTo^{{\pro_{M}}_{2,1}} &
\Imm J^1{\pro_{M}} & \rTo^{{{\pro_{M}}}_{1,0}} & M\times E
  \\
& & \dTo_{\pro_{E}^{(r)}} & & & & \dTo_{\pro_{E}^{(1)}} & &
\dTo_{{\pro_{E}}}
  \\
\ldots & \rTo^{\pi_{r+1,r}} & J^r(E,n) & \rTo^{\pi_{r,r-1}} &
\ldots & \rTo^{\pi_{2,1}} & J^1(E,n) & \rTo^{\pi_{1,0}} & E
\end{diagram}
\end{equation}

\begin{remark}
In \cite{KMS93} the space $J^r(E,n)$ is constructed starting from
$\Imm J^r{\pro_{M}}$ by means of the action of the differential
group $G^r_n$. The previous diagram is an alternative way to get
$J^r(E,n)$ from $\Imm J^r{\pro_{M}}$ by using a generic manifold
$M$ as space of parameters rather than $\R^n$.
\end{remark}

\begin{proposition}\label{prop.comm.diagram}
The following diagram
\begin{equation*}%\label{eq.diagram.comm.for.covering}
\begin{diagram}
  \Imm J^2{\pro_{M}} & \rTo^{\pro_{E}^{(2)}} & J^2(E,n)
  \\
  \uTo^{\ddot{\Gamma}_{{\pro_{M}}}} & & \uTo_{\ddot{\Gamma}}
  \\
  \Imm J^1{\pro_{M}} & \rTo^{\pro_{E}^{(1)}} & J^1(E,n)
\end{diagram}
\end{equation*}
is commutative.
\end{proposition}
\begin{proof}
We calculate the expression of $\pro_{E}^{(2)}$ in a chart where
the determinant of the Jacobian matrix
$\left(u^\xi_{x^\lambda}\right)$ is not zero. We have that
\begin{multline*}%\label{eq.local.exp.of.tau^2}
\pro_{E}^{(2)}(x^\lambda,u^\xi,u^j,u^\xi_{x^\lambda},u^j_{x^\lambda},
u^\xi_{x^\lambda x^\eta},u^j_{x^\lambda x^\eta})=
\\
\left( u^\xi,u^j,x^\lambda_{u^\xi}u^j_{x^\lambda}, u^j_{x^\lambda
x^\mu}x^\lambda_{u^\xi}x^\mu_{u^\eta}-u^j_{x^\lambda}
x^\lambda_{u^\alpha} u^\alpha_{x^\beta x^\gamma} x^\gamma_{u^\xi}
x^\beta_{u^\eta} \right),
\end{multline*}
where $\left(x^\lambda_{u^\xi}\right)$ is the inverse matrix of
$\left(u^\xi_{x^\lambda}\right)$. Then if we take into
consideration diagram \eqref{equ.comm.diagram.tau} and the local
expressions of $\ddot{\Gamma}_{{\pro_{M}}}$ and $\ddot{\Gamma}$,
the proposition follows.
\end{proof}
\begin{theorem}\label{th.covering.equation}
The equation $\cS_{{\pro_{M}}}$ covers the equation $\cS$, with
$\pro_{E}^{(1)}$ as covering map.
\end{theorem}
\begin{proof}
It is sufficient to check that the first prolongation
$\pro_{E}^{(1)}$ sends the distribution
$R\circ\ddot{\Gamma}_{{\pro_{M}}}$ into the distribution
$R\circ\ddot{\Gamma}$. In fact, in view of proposition
\ref{prop.comm.diagram} and of \eqref{equ.prolongations.of.tau},
we have that
\[
R_{\ddot{\Gamma}(\pro_{E}^{(1)}(\tilde{\theta}_1))}=
R_{\pro_{E}^{(2)}\ddot{\Gamma}_{{\pro_{M}}}(\tilde{\theta}_1)}=
T_{\tilde{\theta}_1}\pro_{E}^{(1)}\left(R_{\ddot{\Gamma}_{{\pro_{M}}}(\tilde{\theta}_1)}\right).
\]
\end{proof}

If $R\circ\ddot{\Gamma}$ and $R\circ\ddot{\Gamma}_{{\pro_{M}}}$
are integrable distributions, in view of theorem
\ref{th.formally.integrable}, we have that the previous covering
satisfies the definition given in section
\ref{sec.preliminary.definitions}. In particular, if $n=1$ we find
the results of \cite{MaThesis, MaAnkara}.

\begin{theorem}
The connection $\dot{\Gamma}$ is the quotient connection of
$\dot{\Gamma}_{\pro_{M}}$ via $\pro_{E}^{(1)}$. More precisely the
following diagram
\[
\begin{diagram}
  TJ^1{\pro_{M}} & \rTo^{\upsilon_{\dot{\Gamma}_{\pro_{M}}}} & \ker T{\pro_{M}}_{1,0}
  \\
  \dTo^{T\pro_{E}^{(1)}} & & \dTo_{T\pro_{E}^{(1)}}
  \\
  TJ^1(E,n) & \rTo^{\upsilon_{\dot{\Gamma}}} & \ker T\pi_{1,0}
\end{diagram}
\]
is commutative.
\end{theorem}
\begin{proof}
The local expression of $T\pro_{E}^{(1)}$ is
\[
T\pro_{E}^{(1)}=du^A\otimes\partial_{u^A} -
x^\beta_{u^\xi}x^\lambda_{u^\alpha}u^j_{x^\beta}
du^\xi_{x^\lambda}\otimes\partial_{u^j_\alpha} +
x^\lambda_{u^\alpha}du^i_{x^\lambda}\otimes\partial_{u^i_\alpha}.
\]
The proposition follows in view of
\eqref{eq.local.expr.tangent.and.vertical.valued.forms},
\eqref{eq.local.expr.of.upsilon.dot.nabla} and
\eqref{eq.local.expr.of.dot.nabla.varrho}.
\end{proof}

\section{Example: The case M=$\R^n$}\label{sec.caso.particolare}

Now we analyze the case in which the space of parameters $M$ is
equal to $\R^n$.

\smallskip
The equation $\cS_{\pro_{\R^n}}$ admits always an $n^2+n$
dimensional Lie group of symmetries. More precisely we have the
following
\begin{proposition}
The affine group $\Aff(\R^n)$ of $\R^n$ is a Lie subgroup of the
Lie group of external symmetries of $\cS_{\pro_{\R^n}}$.
\end{proposition}
\begin{proof}
The affine group $\Aff(\R^n)$ induces a Lie group action on
$J^r{\pro_{\R^n}}$. Namely, for any $g\in\Aff(\R^n)$ we can define
the following map:
\begin{equation*}
g\times \id\colon \R^n\times E \rTo  \R^n\times E \,,\quad (x,p)
\mto (g(x),p).
\end{equation*}
Then we can lift it to a Lie transformation $(g\times \id)^{(r)}$
of $J^r{\pro_{\R^n}}$. If $g$ is given locally by
$\tilde{x}^\lambda=a^\lambda_\xi x^\xi+b^\lambda$ with
$a^\lambda_\xi$, $b^{\lambda}\in\R$, we have the following
expression:
\begin{equation}\label{eq.local.expr.of.affine.action}
(g\times\id)^{(2)}(x^\lambda, u^A, u^A_{x^\lambda}, u^A_{x^\lambda
x^\beta})= \left(a^\lambda_\xi
x^\xi+b^\lambda,u^A,u^A_{x^\xi}x^\xi_{\tilde{x}^\lambda},
u^A_{x^\alpha x^\xi} x^\alpha_{\tilde{x}^\lambda}
x^\xi_{\tilde{x}^\beta} \right).
\end{equation}
Now it is a straightforward computation to realize that
$(g\times\id)^{(2)}$ preserves the equation $\cS_{\pro_{\R^n}}$,
and then it is a symmetry.
\end{proof}
Let us note that in the case $n=1$ we recover the affine
parametrization of the geodesic equation.

\begin{proposition}
$\Imm J^1{\pro_{\R^n}}\big/\Aff(\R^n)\simeq J^1(E,n)$.
\end{proposition}
\begin{proof}
Let us denote by $\{[s]^1_{x}\}$ the orbit of $[s]^1_{x}$ with
respect to the action of $\Aff(\R^n)$. Then it follows from
expression \eqref{eq.local.expr.of.affine.action} that the map
\[
\{[s]^1_{x}\} \in \Imm J^1{\pro_{\R^n}}\big/\Aff(\R^n) \mto
\pro_{E}^{(1)}([s]^1_{x})\in J^1(E,n)
\]
is well defined and bijective.
\end{proof}

Finally we have the following

\begin{theorem}
The equation $\cS$ is obtained by factoring the equation
$\cS_{\pro_{\R^n}}$ by the subgroup $\Aff(\R^n)$ of the group of
external symmetries of $\cS_{\pro_{\R^n}}$.
\end{theorem}
\begin{proof}
It is sufficient to consider the results of this section and
theorem \ref{th.covering.equation}.
\end{proof}

\section{Perspectives}
Our next target is to introduce the Weyl projective tensor by
using the approach of this paper, and generalize to the
Grassmannian case. Also, we could apply our results to problems
arising from physics. For instance, we could find contact
projective symmetries, of special form, of the equation of motion
of one relativistic particle.

\paragraph{Acknowledgements.} I would like to thank A.M. Verbovetsky for the continuous support
and V.S. Matveev for stimulating discussions. I would like also to
thank the Universit\`a di Lecce and GNSAGA, which partially
supported this research.

%%%%%%%%%%%%%%%%%%%%%%%%%%%%%%%%%%%%%%%%%%%%%%%%%%%%%%%%%%%%%%%%%%%%%%%%%%%
%%%%%%%%%%%%%%%%%%%%%%%%%%%%%%%%%%%%%%%%%%%%%%%%%%%%%%%%%%%%%%%%%%%%%%%%%%%
%%%%%%%%%%%%%%%%%%%%%%%%%%%%%%%%%%%%%%%%%%%%%%%%%%%%%%%%%%%%%%%%%%%%%%%%%%%


\begin{thebibliography}{AB}
  \small


\bibitem{Aminova01} A. V. Aminova, Projective transformations and
symmetries of differential equations, Sbornik: Mathematics,
\textbf{186}, 1711--1726.

\bibitem{Aminova04} A. V. Aminova, Projective transformations
of pseudo-Riemannian manifolds. Geometry, 9. J. Math. Sci. (N.
Y.), \textbf{113}, no. 3, 2003, 367--470.

%\bibitem{And} I. M. Anderson, The variational bicomplex, book preprint,
%\url{http://www.math.usu.edu/~fg_mp}.

%\bibitem{AVL91} D. V. Alekseevsky, A. M. Vinogradov and V. V.  Lychagin,
%  \textit{Basic ideas and concepts of differential geometry}, Geometry~I.
%  Encycl.\ Math.\ Sci., Vol.\ 28, Springer-Verlag, Berlin, 1991.

\bibitem{Many99} A. V. Bocharov, V. N. Chetverikov, S. V.  Duzhin, N.  G.
  Khor{\cprime}kova, I. S.  Krasil{\cprime}shchik, A.  V. Samokhin, Yu.\ N.
  Torkhov, A. M. Verbovetsky A. M. Vinogradov, \textit{Symmetries and
    Conservation Laws for Differential Equations of Mathematical Physics}, I.
  S. Krasil{\cprime}shchik and A. M.  Vinogradov eds., Amer.\ Math.\ Soc.,
  1999.

\bibitem{Cartan24} E. Cartan, Sur les vari\'et\'e \`a connexions
projective, \textit{Bull. Soc. Math. France}, \textbf{52}, 1924,
205--241.

%\bibitem{Dieu60} J. Dieudonn\'e \textit{Foundation of Modern
%Analysis}, New York and London: Academic Press, 1960.

\bibitem{Dhooghe94} P.F. Dhooghe, Grassmannian structure on manifolds,
\textit{Bull. Belg. Math Soc}, \textbf{1}, 1994, 597-622.

\bibitem{Ehr52} C. Ehresmann, Les prolongements d'une variete differentiable.
  IV. Elements de contact et elements d'enveloppe.  C. R. Acad.\ Sci.\
  Paris,
  \textbf{234}, 1952, 1028--1030.

\bibitem{Eisen22} C. Eisenhart, Spaces with corresponding paths,
\textit{Proc. Nat. Acad. Sci. USA}, \textbf{8}, 1922, 233--238.

%\bibitem{Fat02} L. Fatibene, M. Ferraris, M. Francaviglia, R.G.
%McLenaghan, Generalized symmetries in mechanics and field
%theories, \textit{Journal of Mathematical Physics}, \textbf{43},
%no.6, 2002, 3147-3161.

\bibitem{GrKr98} D. R. Grigore and D. Krupka, Invariant of velocities and
  higher-order Grassmann bundles, \emph{J. Geom.\ Phys.}, \textbf{24}, 1998,
  244--264.

%\bibitem{Hir76} M. H. Hirsch, \textit{Differential Topology},
%Springer-Verlag, 1976.

\bibitem{JanyModyGalilei} A. Jadczyk, M. Modugno, A scheme for Galilei general
 relativistic quantum mechanics, General relativity and gravitational physics
(Bardonecchia, 1992), 319--337, World Sci. Publishing, River Edge, NJ, 1994.

\bibitem{JanyModyBrno} J. Jany\u{s}ka, M. Modugno, Relations between linear
connections on the tangent bundle and connections on the jet
bundle of a fibred manifold, Arch. Math. (Brno), \textbf{32}, 1996, 281--288.

\bibitem{JanyModyRoma} J. Jany\u{s}ka, M. Modugno, An outline of covariant quantum mechanics,
Proceedings of 15th SIGRAV Conference on General Relativity and
Gravitational Physics, (Roma 2002), 281--288.

\bibitem{JanyModyGeneral} J. Jany\u{s}ka, M. Modugno, Classical particle phase
space in general relativity, Proc. Conf. Diff. Geom. and its
Appl., (Brno, 1995), 573--602.

%\bibitem{JanyModyBook} J. Jany\u{s}ka, M. Modugno, Covariant Quantum
%Mechanics, book in preparation, 2005.

\bibitem{Kling82} W. Klingenberg, \emph{Riemannian geometry}, Berlin: Walter de Gruyter, 1982.

%\bibitem{KobTransf} S. Kobayashi: \textit{Transformation group in differential
%geometry}, Springer-Verlag, 1972.

\bibitem{KobNagano} S. Kobayashi and T. Nagano, On projective
connection, \textit{J. Math. Mech}, \textbf{13}, 1964, 215--236.

\bibitem{KMS93} I. Kol\'a\v r, P. Michor and J. Slov\'ak, \textit{Natural
    Operations in Differential Geometry}, Springer-Verlag, 1993.

%\bibitem{Kol98} I. Kol\'{a}\v{r}, Affine structure on Weil bundles,
%  \textit{Nagoya Math.\ J.}, \textbf{158}, 2000, 99-106.

%\bibitem{KLV86} I. S. Krasil{\cprime}shchik, V. V. Lychagin and A. M.
%  Vinogradov, \emph{Geometry of Jet Spaces and Non--linear Partial Differential
%    Equations}, Gordon and Breach, New York, 1986.

\bibitem{KrVe98} I. S. Krasil{\cprime}shchik and A. M. Verbovetsky,
  \emph{Homological methods in equations of mathematical physics}, Open
  Education and Sciences, Opava (Czech Rep.), 1998 math.DG/9808130.

\bibitem{KrVi84} I. S. Krasil{\cprime}shchik and A. M. Vinogradov,
  Non local symmetries and the theory of covering:
  An addendum to A.M. Vinogradov's `Local symmetries and conservation laws',
  Acta Appl. Math, \textbf{2}, 1984, 79-96.

\bibitem{KrVi89} I. S. Krasil{\cprime}shchik and A. M. Vinogradov,
  Nonlocal trends in the geometry of differential equations:
  symmetries, conservation laws and B\"{a}cklund transformation,
  Acta Appl. Math \textbf{15}, 1989, 161-209.

%\bibitem{Krey59} E. Kreyszig, \emph{Differential geometry}, University of Toronto Press, 1959.

\bibitem{Kru01} D. Krupka, Global variational functionals in fibered spaces,
  \emph{Nonlinear analysis}, \textbf{47}, 2001, 2633--2642.

\bibitem{MaThesis} G. Manno, \textit{Jet methods for the finite order
variational sequence and the geodesic equation}, Ph.D. thesis,
September 2003.

\bibitem{MaAnkara} G. Manno, \textit{The geometry of the geodesic equation in the framework of jets of
submanifolds}, Conference proceedings of AIP, \textbf{729}, 2004.

\bibitem{MaVi02} G. Manno and R. Vitolo, Variational sequences on finite order
  jets of submanifolds, Proc. Conf. Diff. Geom. and its
Appl., (Opava 2001), 435--446.

\bibitem{MaVi03} G. Manno and R. Vitolo, Some cohomological aspects of the calculus of
variations on finite-order jets of submanifolds, submitted to
\emph{Duke Math. Journal}.

\bibitem{MaVi04} G. Manno and R. Vitolo, Relativistic mechanics, cosymplectic manifolds and
      symmetries, Note di Matematica, \textbf{23}, no. 2, 2004.

%\bibitem{MaVi05} G. Manno and R. Vitolo, Contact
%projective symmetries of the motion of one relativistic particle,
%in preparation.

\bibitem{MoVi94} M. Modugno, A. M. Vinogradov, Some Variations on the Notion of
  Connection, \emph{Ann.\ di Mat.\ Pura ed Appl.} IV, Vol.\ CLXVII, 1994,
  33--71.

\bibitem{Olver01} P. J. Olver, \textit{Application of Lie groups
to differential equations. Second edition}, Springer-Verlag, New
York, 1993.

%\bibitem{Osb82} H. Osborn, \emph{Vector Bundles}, New York ; London : Academic Press, 1982.

%\bibitem{Pom78} J. F. Pommaret, \textit{Systems of partial differential equations and Lie pseudogroups},
%Gordon and Breach Science Publishers, New York, 1978.

\bibitem{Sau89} D. J. Saunders, \textit{The Geometry of Jet Bundles}, Cambridge
  Univ.\ Press, 1989.

\bibitem{Thomas25} T.Y. Thomas, On the projective and
equi-projective geometry of path, \textit{Proc. Nat. Acad. Sci.},
\textbf{11}, 1925, 199--203.

\bibitem{Thomas26} T.Y. Thomas, A projective theory of affinely connected manifolds,
\textit{Math. Zeit.}, \textbf{25}, 1926, 723--733.

\bibitem{Veblen22} O. Veblen, Projective and affine geometry of
paths, \textit{Proc. Nat. Acad. Sci. USA}, \textbf{8}, 1922,
347--350.

\bibitem{Vin81} A. M. Vinogradov, Local symmetries and conservation laws,
 \emph{Acta Applicandae Mathematicae},  \textbf{2}, 1981, 21--78.

\bibitem{Vin84} A. M. Vinogradov, The $\cal{C}$-spectral Sequence, Lagrangian
  Formalism and Conservation Laws I and II, \emph{Jour.\ of Math.\ Analysis and
    Appl.}, \textbf{100}, no. 1, 1984, 1--129.

\bibitem{Vin88} A. M. Vinogradov, An informal introduction to the geometry of
  jet spaces, \emph{Rend.\ Se\-mi\-na\-ri Fac.\ Sci.\ Univ.\
  Ca\-glia\-ri}, \textbf{58}, 1988, 301--333.

\bibitem{VinII} A. M. Vinogradov, \textit{Cohomological Analysis of Partial
Differential Equations and Secondary Calculus}, Amer.\ Math.\
Soc., 2001.

\bibitem{Weyl21} H. Weyl, Zur Infinitesimalgeometrie; Einordnung der
projektiven und der knoformen Auffassung, \textit{Gottingen
Nachr.}, 1921, 99--122.

%\bibitem{Yano70} K. Yano, \textit{Integral Formulas in Riemannian
%Geometry}, New York, 1970.

\end{thebibliography}
\end{document}